# ON THE "DEGREES OF FREEDOM" OF THE LASSO

By Hui Zou, Trevor Hastie and Robert Tibshirani

*University of Minnesota, Stanford University and Stanford University*

We study the effective degrees of freedom of the lasso in the framework of Stein's unbiased risk estimation (SURE). We show that the number of nonzero coefficients is an unbiased estimate for the degrees of freedom of the lasso—a conclusion that requires no special assumption on the predictors. In addition, the unbiased estimator is shown to be asymptotically consistent. With these results on hand, various model selection criteria—$C_p$, AIC and BIC—are available, which, along with the LARS algorithm, provide a principled and efficient approach to obtaining the optimal lasso fit with the computational effort of a single ordinary least-squares fit.

**1. Introduction.** The lasso is a popular model building technique that simultaneously produces accurate and parsimonious models (Tibshirani [22]). Suppose $\mathbf{y} = (y_1, \ldots, y_n)^T$ is the response vector and $\mathbf{x}_j = (x_{1j}, \ldots, x_{nj})^T$, $j = 1, \ldots, p$, are the linearly independent predictors. Let $\mathbf{X} = [\mathbf{x}_1, \ldots, \mathbf{x}_p]$ be the predictor matrix. Assume the data are standardized. The lasso estimates for the coefficients of a linear model are obtained by

$$(1.1) \qquad \hat{\beta} = \arg\min_{\beta} \left\| \mathbf{y} - \sum_{j=1}^{p} \mathbf{x}_j \beta_j \right\|^2 + \lambda \sum_{j=1}^{p} |\beta_j|,$$

where $\lambda$ is called the lasso regularization parameter. What we show in this paper is that the number of nonzero components of $\hat{\beta}$ is an exact unbiased estimate of the degrees of freedom of the lasso, and this result can be used to construct adaptive model selection criteria for efficiently selecting the optimal lasso fit.

Degrees of freedom is a familiar phrase for many statisticians. In linear regression the degrees of freedom is the number of estimated predictors. Degrees of freedom is often used to quantify the model complexity of a statistical modeling procedure (Hastie and Tibshirani [10]). However, generally speaking, there is no exact correspondence between the degrees of freedom and the number of parameters in the model (Ye [24]). For example, suppose we first find $x_{j^*}$ such that $|\text{cor}(x_{j^*}, y)|$ is the largest among all $x_j, j = 1, 2, \ldots, p$. We then use $x_{j^*}$ to fit a simple linear regression model to predict $y$. There is one parameter in the fitted model, but the degrees of freedom is greater than one, because we have to take into account the stochastic search of $x_{j^*}$.

Stein's unbiased risk estimation (SURE) theory (Stein [21]) gives a rigorous definition of the degrees of freedom for any fitting procedure. Given a model fitting method $\delta$, let $\hat{\boldsymbol{\mu}} = \delta(\mathbf{y})$ represent its fit. We assume that given the $\mathbf{x}$'s, $\mathbf{y}$ is generated according to $\mathbf{y} \sim (\boldsymbol{\mu}, \sigma^2 \mathbf{I})$, where $\boldsymbol{\mu}$ is the true mean vector and $\sigma^2$ is the common variance. It is shown (Efron [4]) that the degrees of freedom of $\delta$ is

$$(1.2) \qquad df(\hat{\boldsymbol{\mu}}) = \sum_{i=1}^{n} \text{cov}(\hat{\mu}_i, y_i)/\sigma^2.$$

For example, if $\delta$ is a linear smoother, that is, $\hat{\boldsymbol{\mu}} = \mathbf{S}\mathbf{y}$ for some matrix $\mathbf{S}$ independent of $\mathbf{y}$, then we have $\text{cov}(\hat{\boldsymbol{\mu}}, \mathbf{y}) = \sigma^2 \mathbf{S}$, $df(\hat{\boldsymbol{\mu}}) = \text{tr}(\mathbf{S})$. SURE theory also reveals the statistical importance of the degrees of freedom. With $df$ defined in (1.2), we can employ the covariance penalty method to construct a $C_p$-type statistic as

$$(1.3) \qquad C_p(\hat{\boldsymbol{\mu}}) = \frac{\|\mathbf{y} - \hat{\boldsymbol{\mu}}\|^2}{n} + \frac{2df(\hat{\boldsymbol{\mu}})}{n}\sigma^2.$$

Efron [4] showed that $C_p$ is an unbiased estimator of the true prediction error, and in some settings it offers substantially better accuracy than cross-validation and related nonparametric methods. Thus degrees of freedom plays an important role in model assessment and selection. Donoho and Johnstone [3] used the SURE theory to derive the degrees of freedom of soft thresholding and showed that it leads to an adaptive wavelet shrinkage procedure called *SureShrink*. Ye [24] and Shen and Ye [20] showed that the degrees of freedom can capture the inherent uncertainty in modeling and frequentist model selection. Shen and Ye [20] and Shen, Huang and Ye [19] further proved that the degrees of freedom provides an adaptive model selection criterion that performs better than the fixed-penalty model selection criteria.

The lasso is a regularization method which does automatic variable selection. As shown in Figure 1 (the left panel), the lasso continuously shrinks the coefficients toward zero as $\lambda$ increases; and some coefficients are shrunk to exactly zero if $\lambda$ is sufficiently large. Continuous shrinkage also often improves the prediction accuracy due to the bias–variance trade-off. Detailed discussions on variable selection via penalization are given in Fan and Li [6], Fan and Peng [8] and Fan and Li [7]. In recent years the lasso has attracted a lot of attention in both the statistics and machine learning communities. It is of great interest to know the degrees of freedom of the lasso for any given regularization parameter $\lambda$ for selecting the optimal lasso model. However, it is difficult to derive the analytical expression of the degrees of freedom of many nonlinear modeling procedures, including the lasso. To overcome the analytical difficulty, Ye [24] and Shen and Ye [20] proposed using a data-perturbation technique to numerically compute an (approximately) unbiased estimate for $df(\hat{\boldsymbol{\mu}})$ when the analytical form of $\hat{\boldsymbol{\mu}}$ is unavailable. The bootstrap (Efron [4]) can also be used to obtain an (approximately) unbiased estimator of the degrees of freedom. This kind of approach, however, can be computationally expensive. It is an interesting problem of both theoretical and practical importance to derive rigorous analytical results on the degrees of freedom of the lasso.

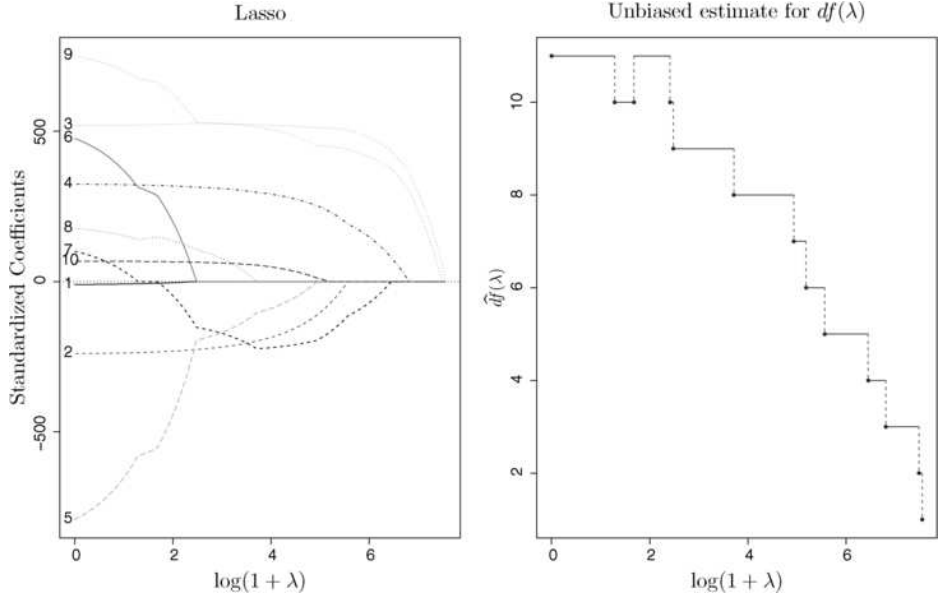

FIG. 1. *Diabetes data with ten predictors. The left panel shows the lasso coefficient estimates $\hat{\beta}_j, j = 1, 2, \ldots, 10$, for the diabetes study. The lasso coefficient estimates are piece-wise linear functions of $\lambda$ (Osborne, Presnell and Turlach [15] and Efron, Hastie, Johnstone and Tibshirani [5]), hence they are piece-wise nonlinear as functions of $\log(1 + \lambda)$. The right panel shows the curve of the proposed unbiased estimate for the degrees of freedom of the lasso.*

In this work we study the degrees of freedom of the lasso in the framework of SURE. We show that for any given $\lambda$ the number of nonzero predictors in the model is an unbiased estimate for the degrees of freedom. This is a finite-sample exact result and the result holds as long as the predictor matrix is a full rank matrix. The importance of the exact finite-sample unbiasedness is emphasized in Efron [4], Shen and Ye [20] and Shen and Huang [18]. We show that the unbiased estimator is also consistent. As an illustration, the right panel in Figure 1 displays the unbiased estimate for the degrees of freedom as a function of $\lambda$ for the diabetes data (with ten predictors).

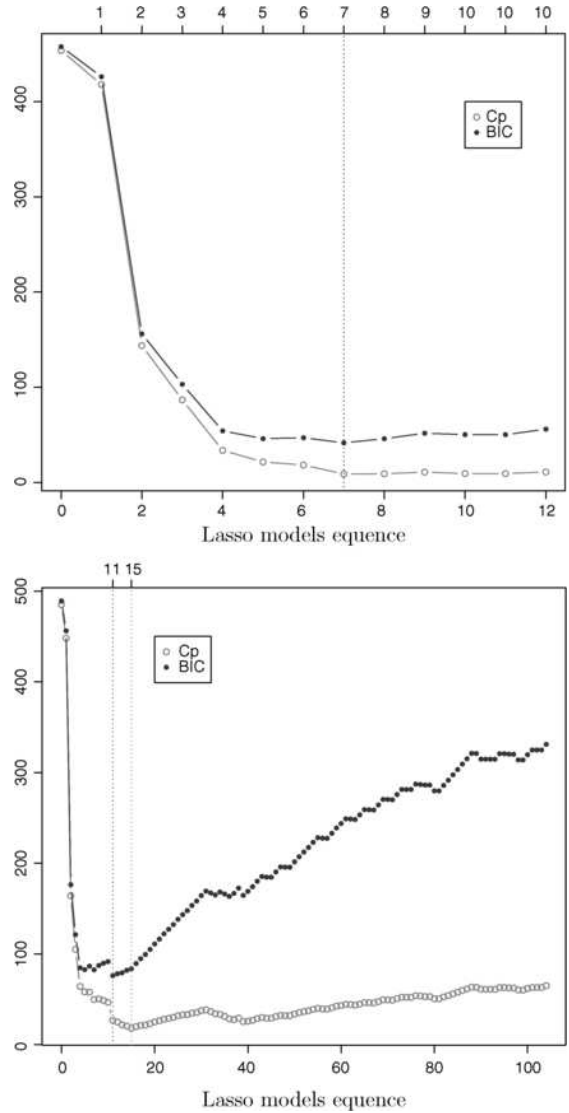

FIG. 2. *The diabetes data: $C_p$ and BIC curves with ten (top) and 64 (bottom) predictors. In the top panel $C_p$ and BIC select the same model with seven nonzero coefficients. In the bottom panel, $C_p$ selects a model with 15 nonzero coefficients and BIC selects a model with 11 nonzero coefficients.*

The unbiased estimate of the degrees of freedom can be used to construct $C_p$ and BIC type model selection criteria. The $C_p$ (or BIC) curve is easily obtained once the lasso solution paths are computed by the LARS algorithm (Efron, Hastie, Johnstone and Tibshirani [5]). Therefore, with the computational effort of a single OLS fit, we are able to find the optimal lasso fit using our theoretical results. Note that $C_p$ is a finite-sample result and relies on its unbiasedness for prediction error as a basis for model selection (Shen and Ye [20], Efron [4]). For this purpose, an unbiased estimate of the degrees of freedom is sufficient. We illustrate the use of $C_p$ and BIC on the diabetes data in Figure 2, where the selected models are indicated by the broken vertical lines.

The rest of the paper is organized as follows. We present the main results in Section 2. We construct model selection criteria—$C_p$ or BIC—using the degrees of freedom. In Section 3 we discuss the conjecture raised in [5]. Section 4 contains some technical proofs. Discussion is in Section 5.

**2. Main results.** We first define some notation. Let $\hat{\boldsymbol{\mu}}_\lambda$ be the lasso fit using the representation (1.1). $\hat{\mu}_i$ is the $i$th component of $\hat{\boldsymbol{\mu}}$. For convenience, we let $df(\lambda)$ stand for $df(\hat{\boldsymbol{\mu}}_\lambda)$, the degrees of freedom of the lasso. Suppose **M** is a matrix with $p$ columns. Let $\mathcal{S}$ be a subset of the indices $\{1, 2, \ldots, p\}$. Denote by $\mathbf{M}_\mathcal{S}$ the submatrix $\mathbf{M}_\mathcal{S} = [\cdots M_j \cdots]_{j \in \mathcal{S}}$, where $M_j$ is the $j$th column of **M**. Similarly, define $\beta_\mathcal{S} = (\cdots \beta_j \cdots)_{j \in \mathcal{S}}$ for any vector $\beta$ of length $p$. Let $\text{Sgn}(\cdot)$ be the sign function: $\text{Sgn}(x) = 1$ if $x > 0$; $\text{Sgn}(x) = 0$ if $x = 0$; $\text{Sgn}(x) = -1$ if $x = -1$. Let $\mathcal{B} = \{j : \text{Sgn}(\beta)_j \neq 0\}$ be the *active set* of $\beta$, where $\text{Sgn}(\beta)$ is the sign vector of $\beta$ given by $\text{Sgn}(\beta)_j = \text{Sgn}(\beta_j)$. We denote the active set of $\hat{\beta}(\lambda)$ as $\mathcal{B}(\lambda)$ and the corresponding sign vector $\text{Sgn}(\hat{\beta}(\lambda))$ as $\text{Sgn}(\lambda)$. We do not distinguish between the index of a predictor and the predictor itself.

2.1. *The unbiased estimator of* $df(\lambda)$. Before delving into the technical details, let us review some characteristics of the lasso solution (Efron et al. [5]). For a given response vector **y**, there is a *finite* sequence of $\lambda$'s,

$$\lambda_0 > \lambda_1 > \lambda_2 > \cdots > \lambda_K = 0, \tag{2.1}$$

such that:

- For all $\lambda > \lambda_0$, $\hat{\beta}(\lambda) = 0$.
- In the interior of the interval $(\lambda_{m+1}, \lambda_m)$, the active set $\mathcal{B}(\lambda)$ and the sign vector $\text{Sgn}(\lambda)_{\mathcal{B}(\lambda)}$ are constant with respect to $\lambda$. Thus we write them as $\mathcal{B}_m$ and $\text{Sgn}_m$ for convenience.

The active set changes at each $\lambda_m$. When $\lambda$ decreases from $\lambda = \lambda_m - 0$, some predictors with zero coefficients at $\lambda_m$ are about to have nonzero coefficients; thus they join the active set $\mathcal{B}_m$. However, as $\lambda$ approaches $\lambda_{m+1} + 0$ there are possibly some predictors in $\mathcal{B}_m$ whose coefficients reach zero. Hence we call $\{\lambda_m\}$ the *transition points*. Any $\lambda \in [0, \infty) \setminus \{\lambda_m\}$ is called a nontransition point.

THEOREM 1.   $\forall \lambda$ *the lasso fit $\hat{\boldsymbol{\mu}}_\lambda(\mathbf{y})$ is a uniformly Lipschitz function on* $\mathbf{y}$. *The degrees of freedom of $\hat{\boldsymbol{\mu}}_\lambda(\mathbf{y})$ equal the expectation of the effective set* $\mathcal{B}_\lambda$, *that is,*

$$df(\lambda) = E|\mathcal{B}_\lambda|. \qquad (2.2)$$

*The identity (2.2) holds as long as* $\mathbf{X}$ *is full rank, that is,* $\text{rank}(\mathbf{X}) = p$.

Theorem 1 shows that $\widehat{df}(\lambda) = |\mathcal{B}_\lambda|$ is an unbiased estimate for $df(\lambda)$. Thus $\widehat{df}(\lambda)$ suffices to provide an exact unbiased estimate to the true prediction risk of the lasso. The importance of the exact finite-sample unbiasedness is emphasized in Efron [4], Shen and Ye [20] and Shen and Huang [18]. Our result is also computationally friendly. Given any data set, the entire solution paths of the lasso are computed by the LARS algorithm (Efron et al. [5]); then the unbiased estimator $\widehat{df}(\lambda) = |\mathcal{B}_\lambda|$ is easily obtained without any extra effort.

To prove Theorem 1 we shall proceed by proving a series of lemmas whose proofs are relegated to Section 4 for the sake of presentation.

LEMMA 1.   *Suppose $\lambda \in (\lambda_{m+1}, \lambda_m)$. $\hat{\beta}(\lambda)$ are the lasso coefficient estimates. Then we have*

$$\hat{\beta}(\lambda)_{\mathcal{B}_m} = (\mathbf{X}_{\mathcal{B}_m}^T \mathbf{X}_{\mathcal{B}_m})^{-1} \left( \mathbf{X}_{\mathcal{B}_m}^T \mathbf{y} - \frac{\lambda}{2} \text{Sgn}_m \right). \qquad (2.3)$$

LEMMA 2.   *Consider the transition points $\lambda_m$ and $\lambda_{m+1}$, $\lambda_{m+1} \geq 0$. $\mathcal{B}_m$ is the active set in $(\lambda_{m+1}, \lambda_m)$. Suppose $i_{\text{add}}$ is an index added into $\mathcal{B}_m$ at $\lambda_m$ and its index in $\mathcal{B}_m$ is $i^*$, that is, $i_{\text{add}} = (\mathcal{B}_m)_{i^*}$. Denote by $(a)_k$ the $k$th element of the vector $a$. We can express the transition point $\lambda_m$ as*

$$\lambda_m = \frac{2((\mathbf{X}_{\mathcal{B}_m}^T \mathbf{X}_{\mathcal{B}_m})^{-1} \mathbf{X}_{\mathcal{B}_m}^T \mathbf{y})_{i^*}}{((\mathbf{X}_{\mathcal{B}_m}^T \mathbf{X}_{\mathcal{B}_m})^{-1} \text{Sgn}_m)_{i^*}}. \qquad (2.4)$$

*Moreover, if $j_{\text{drop}}$ is a dropped (if there is any) index at $\lambda_{m+1}$ and $j_{\text{drop}} = (\mathcal{B}_m)_{j^*}$, then $\lambda_{m+1}$ can be written as*

$$\lambda_{m+1} = \frac{2((\mathbf{X}_{\mathcal{B}_m}^T \mathbf{X}_{\mathcal{B}_m})^{-1} \mathbf{X}_{\mathcal{B}_m}^T \mathbf{y})_{j^*}}{((\mathbf{X}_{\mathcal{B}_m}^T \mathbf{X}_{\mathcal{B}_m})^{-1} \text{Sgn}_m)_{j^*}}. \qquad (2.5)$$

LEMMA 3. $\forall \lambda > 0$, $\exists$ a null set $\mathcal{N}_\lambda$ which is a finite collection of hyperplanes in $\mathbb{R}^n$. Let $\mathcal{G}_\lambda = \mathbb{R}^n \setminus \mathcal{N}_\lambda$. Then $\forall \mathbf{y} \in \mathcal{G}_\lambda$, $\lambda$ is not any of the transition points, that is, $\lambda \notin \{\lambda(\mathbf{y})_m\}$.

LEMMA 4. $\forall \lambda$, $\hat{\beta}_\lambda(\mathbf{y})$ is a continuous function of $\mathbf{y}$.

LEMMA 5. Fix any $\lambda > 0$ and consider $\mathbf{y} \in \mathcal{G}_\lambda$ as defined in Lemma 3. The active set $\mathcal{B}(\lambda)$ and the sign vector $\mathrm{Sgn}(\lambda)$ are locally constant with respect to $\mathbf{y}$.

LEMMA 6. Let $\mathcal{G}_0 = \mathbb{R}^n$. Fix an arbitrary $\lambda \geq 0$. On the set $\mathcal{G}_\lambda$ with full measure as defined in Lemma 3, the lasso fit $\hat{\boldsymbol{\mu}}_\lambda(\mathbf{y})$ is uniformly Lipschitz. Precisely,

(2.6) $\quad \|\hat{\boldsymbol{\mu}}_\lambda(\mathbf{y} + \Delta \mathbf{y}) - \hat{\boldsymbol{\mu}}_\lambda(\mathbf{y})\| \leq \|\Delta \mathbf{y}\| \quad$ for sufficiently small $\Delta \mathbf{y}$.

Moreover, we have the divergence formula

(2.7) $$\nabla \cdot \hat{\boldsymbol{\mu}}_\lambda(\mathbf{y}) = |\mathcal{B}_\lambda|.$$

PROOF OF THEOREM 1. Theorem 1 is obviously true for $\lambda = 0$. We only need to consider $\lambda > 0$. By Lemma 6 $\hat{\boldsymbol{\mu}}_\lambda(\mathbf{y})$ is uniformly Lipschitz on $\mathcal{G}_\lambda$. Moreover, $\hat{\boldsymbol{\mu}}_\lambda(\mathbf{y})$ is a continuous function of $\mathbf{y}$, and thus $\hat{\boldsymbol{\mu}}_\lambda(\mathbf{y})$ is uniformly Lipschitz on $\mathbb{R}^n$. Hence $\hat{\boldsymbol{\mu}}_\lambda(\mathbf{y})$ is *almost differentiable*; see Meyer and Woodroofe [14] and Efron et al. [5]. Then (2.2) is obtained by invoking Stein's lemma (Stein [21]) and the divergence formula (2.7). □

2.2. *Consistency of the unbiased estimator $\widehat{df}(\lambda)$.* In this section we show that the obtained unbiased estimator $\widehat{df}(\lambda)$ is also consistent. We adopt the similar setup in Knight and Fu [12] for the asymptotic analysis. Assume the following two conditions:

1. $y_i = \mathbf{x}_i \beta^* + \varepsilon_i$, where $\varepsilon_1, \ldots, \varepsilon_n$ are i.i.d. normal random variables with mean 0 and variance $\sigma^2$, and $\beta^*$ denotes the fixed unknown regression coefficients.
2. $\frac{1}{n} \mathbf{X}^T \mathbf{X} \to \mathbf{C}$, where $\mathbf{C}$ is a positive definite matrix.

We consider minimizing an objective function $Z_\lambda(\beta)$ defined as

(2.8) $$Z_\lambda(\beta) = (\beta - \beta^*)^T \mathbf{C} (\beta - \beta^*) + \lambda \sum_{j=1}^p |\beta_j|.$$

Optimizing (2.8) is a lasso type problem: minimizing a quadratic objective function with an $\ell_1$ penalty. There are also a *finite* sequence of transition points $\{\lambda_{*m}\}$ associated with optimizing (2.8).

THEOREM 2. *If $\frac{\lambda_n^*}{n} \to \lambda^* > 0$, where $\lambda^*$ is a nontransition point such that $\lambda^* \neq \lambda_{*m}$ for all $m$, then $\widehat{df}(\lambda_n^*) - df(\lambda_n^*) \to 0$ in probability.*

PROOF OF THEOREM 2. Consider $\hat{\beta}^* = \arg\min_\beta Z_{\lambda^*}(\beta)$ and let $\hat{\beta}^{(n)}$ be the lasso solution given in (1.1) with $\lambda = \lambda_n^*$. Denote $\mathcal{B}^{(n)} = \{j : \hat{\beta}_j^{(n)} \neq 0, 1 \leq j \leq p\}$ and $\mathcal{B}^* = \{j : \hat{\beta}_j^* \neq 0, 1 \leq j \leq p\}$. We want to show $P(\mathcal{B}^{(n)} = \mathcal{B}^*) \to 1$. First, let us consider any $j \in \mathcal{B}^*$. By Theorem 1 in Knight and Fu [12] we know that $\hat{\beta}^{(n)} \to_p \hat{\beta}^*$. Then the continuous mapping theorem implies that $\mathrm{Sgn}(\hat{\beta}_j^{(n)}) \to_p \mathrm{Sgn}(\hat{\beta}_j^*) \neq 0$, since $\mathrm{Sgn}(x)$ is continuous at all $x$ but zero. Thus $P(\mathcal{B}^{(n)} \supseteq \mathcal{B}^*) \to 1$. Second, consider any $j' \notin \mathcal{B}^*$. Then $\hat{\beta}_{j'}^* = 0$. Since $\hat{\beta}^*$ is the minimizer of $Z_{\lambda^*}(\beta)$ and $\lambda^*$ is not a transition point, by the Karush–Kuhn–Tucker (KKT) optimality condition (Efron et al. [5], Osborne, Presnell and Turlach [15]), we must have

$$(2.9) \qquad \lambda^* > 2|\mathbf{C}_{j'}(\beta^* - \hat{\beta}^*)|,$$

where $\mathbf{C}_{j'}$ is the $j'$th row vector of $\mathbf{C}$. Let $r^* = \lambda^* - 2|\mathbf{C}_{j'}(\beta^* - \hat{\beta}^*)| > 0$. Now let us consider $r_n = \lambda_n^* - 2|\mathbf{x}_{j'}^T(\mathbf{y} - \mathbf{X}\hat{\beta}_n^*)|$. Note that

$$(2.10) \qquad \mathbf{x}_{j'}^T(\mathbf{y} - \mathbf{X}\hat{\beta}_n^*) = \mathbf{x}_{j'}^T\mathbf{X}(\beta^* - \hat{\beta}_n^*) + \mathbf{x}_{j'}^T\varepsilon.$$

Thus $\frac{r_n^*}{n} = \frac{\lambda_n^*}{n} - 2|\frac{1}{n}\mathbf{x}_{j'}^T\mathbf{X}(\beta^* - \hat{\beta}_n^*) + \mathbf{x}_{j'}^T\varepsilon/n|$. Because $\hat{\beta}^{(n)} \to_p \hat{\beta}^*$ and $\mathbf{x}_{j'}^T\varepsilon/n \to_p 0$, we conclude $\frac{r_n^*}{n} \to_p r^* > 0$. By the KKT optimality condition, $r_n^* > 0$ implies $\hat{\beta}_{j'}^{(n)} = 0$. Thus $P(\mathcal{B}^* \supseteq \mathcal{B}^n) \to 1$. Therefore $P(\mathcal{B}^{(n)} = \mathcal{B}^*) \to 1$. Immediately we see $\widehat{df}(\lambda_n^*) \to_p |\mathcal{B}^*|$. Then invoking the dominated convergence theorem we have

$$(2.11) \qquad df(\lambda_n^*) = E[\widehat{df}(\lambda_n^*)] \to |\mathcal{B}^*|.$$

So $\widehat{df}(\lambda_n^*) - df(\lambda_n^*) \to_p 0$. $\square$

2.3. *Numerical experiments.* In this section we check the validity of our arguments by a simulation study. Here is the outline of the simulation. We take the 64 predictors in the diabetes data set, which include the quadratic terms and interactions of the original ten predictors. The positive cone condition is violated on the 64 predictors (Efron et al. [5]). The response vector $\mathbf{y}$ is used to fit an OLS model. We compute the OLS estimates $\hat{\beta}_{\mathrm{ols}}$ and $\hat{\sigma}_{\mathrm{ols}}^2$. Then we consider a synthetic model,

$$(2.12) \qquad \mathbf{y}^* = \mathbf{X}\beta + N(0,1)\sigma,$$

where $\beta = \hat{\beta}_{\mathrm{ols}}$ and $\sigma = \hat{\sigma}_{\mathrm{ols}}$.

Given the synthetic model, the degrees of freedom of the lasso can be numerically evaluated by Monte Carlo methods. For $b = 1, 2, \ldots, B$, we independently simulate $\mathbf{y}^*(b)$ from (2.12). For a given $\lambda$, by the definition of $df$, we need to evaluate $\text{cov}_i = \text{cov}(\hat{\mu}_i, y_i^*)$. Then $df = \sum_{i=1}^n \text{cov}_i / \sigma^2$. Since $E[y_i^*] = (\mathbf{X}\beta)_i$ and note that $\text{cov}_i = E[(\hat{\mu}_i - a_i)(y_i^* - (\mathbf{X}\beta)_i)]$ for any fixed known constant $a_i$. Then we compute

$$\widehat{\text{cov}}_i = \frac{\sum_{b=1}^B (\hat{\mu}_i(b) - a_i)(y_i^*(b) - (\mathbf{X}\beta)_i)}{B} \tag{2.13}$$

and $df = \sum_{i=1}^n \widehat{\text{cov}}_i / \sigma^2$. Typically $a_i = 0$ is used in Monte Carlo calculation. In this work we use $a_i = (\mathbf{X}\beta)_i$, for it gives a Monte Carlo estimate for $df$ with smaller variance than that given by $a_i = 0$. On the other hand, we evaluate $E|\mathcal{B}_\lambda|$ by $\sum_{b=1}^B \widehat{df}(\lambda)_b / B$. We are interested in $E|\mathcal{B}_\lambda| - df(\lambda)$. Standard errors are calculated based on the $B$ replications. Figure 3 shows very convincing pictures to support the identity (2.2).

2.4. *Adaptive model selection criteria.* The exact value of $df(\lambda)$ depends on the underlying model according to Theorem 1. It remains unknown to us unless we know the underlying model. Our theory provides a convenient unbiased and consistent estimate of the unknown $df(\lambda)$. In the spirit of SURE theory, the good unbiased estimate for $df(\lambda)$ suffices to provide an unbiased estimate for the prediction error of $\hat{\boldsymbol{\mu}}_\lambda$ as

$$C_p(\hat{\boldsymbol{\mu}}) = \frac{\|\mathbf{y} - \hat{\boldsymbol{\mu}}\|^2}{n} + \frac{2}{n} \widehat{df}(\hat{\boldsymbol{\mu}}) \sigma^2. \tag{2.14}$$

Consider the $C_p$ curve as a function of the regularization parameter $\lambda$. We find the optimal $\lambda$ that minimizes $C_p$. As shown in Shen and Ye [20], this model selection approach leads to an adaptively optimal model which essentially achieves the optimal prediction risk as if the ideal tuning parameter were given in advance.

By the connection between Mallows' $C_p$ (Mallows [13]) and AIC (Akaike [1]), we use the (generalized) $C_p$ formula (2.14) to equivalently define AIC for the lasso,

$$\text{AIC}(\hat{\boldsymbol{\mu}}) = \frac{\|\mathbf{y} - \hat{\boldsymbol{\mu}}\|^2}{n\sigma^2} + \frac{2}{n} \widehat{df}(\hat{\boldsymbol{\mu}}). \tag{2.15}$$

The model selection results are identical by $C_p$ and AIC. Following the usual definition of BIC [16], we propose BIC for the lasso as

$$\text{BIC}(\hat{\boldsymbol{\mu}}) = \frac{\|\mathbf{y} - \hat{\boldsymbol{\mu}}\|^2}{n\sigma^2} + \frac{\log(n)}{n} \widehat{df}(\hat{\boldsymbol{\mu}}). \tag{2.16}$$

AIC and BIC possess different asymptotic optimality. It is well known that AIC tends to select the model with the optimal prediction performance,

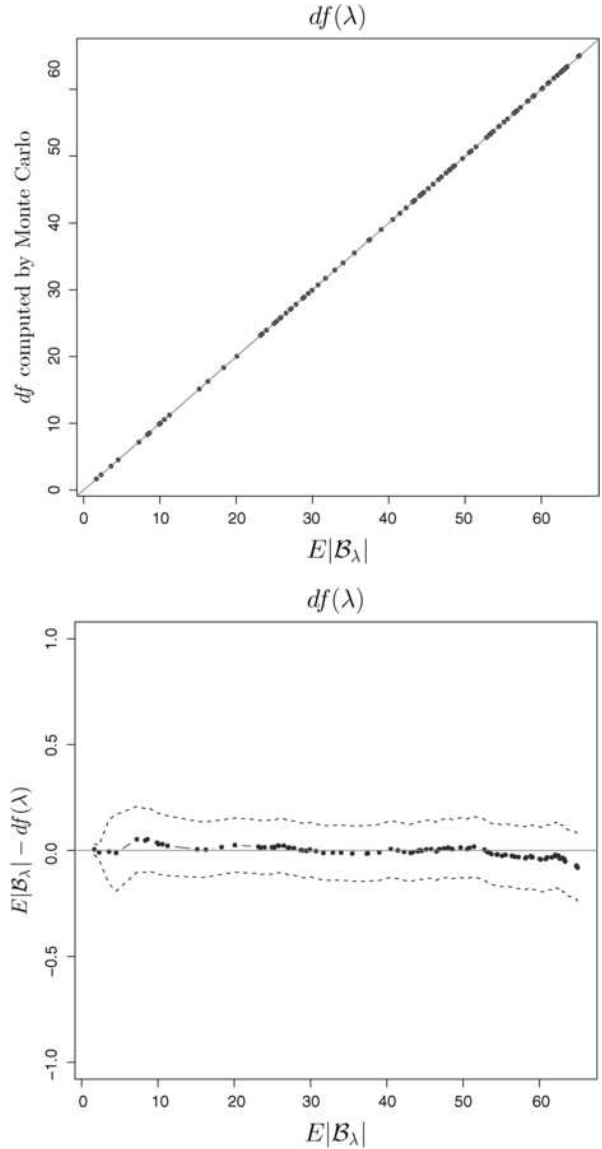

Fig. 3. *The synthetic model with the 64 predictors in the diabetes data. In the top panel we compare $E|\mathcal{B}_\lambda|$ with the true degrees of freedom $df(\lambda)$ based on $B = 20000$ Monte Carlo simulations. The solid line is the $45°$ line (the perfect match). The bottom panel shows the estimation bias and its point-wise $95\%$ confidence intervals are indicated by the thin dashed lines. Note that the zero horizontal line is well inside the confidence intervals.*

while BIC tends to identify the true sparse model if the true model is in the candidate list; see Shao [17], Yang [23] and references therein. We suggest

DEGREES OF FREEDOM OF THE LASSO    11using BIC as the model selection criterion when the sparsity of the model is our primary concern.

Using either AIC or BIC to find the optimal lasso model, we are facing an optimization problem,

$$\lambda(\text{optimal}) = \arg\min_{\lambda} \frac{\|\mathbf{y} - \hat{\boldsymbol{\mu}}_\lambda\|^2}{n\sigma^2} + \frac{w_n}{n}\widehat{df}(\lambda), \tag{2.17}$$

where $w_n = 2$ for AIC and $w_n = \log(n)$ for BIC. Since the LARS algorithm efficiently solves the lasso solution for all $\lambda$, finding $\lambda(\text{optimal})$ is attainable in principle. In fact, we show that $\lambda(\text{optimal})$ is one of the transition points, which further facilitates the searching procedure.

THEOREM 3. *To find $\lambda(\text{optimal})$, we only need to solve*

$$m^* = \arg\min_{m} \frac{\|\mathbf{y} - \hat{\boldsymbol{\mu}}_{\lambda_m}\|^2}{n\sigma^2} + \frac{w_n}{n}\widehat{df}(\lambda_m); \tag{2.18}$$

*then $\lambda(\text{optimal}) = \lambda_{m^*}$.*

PROOF.  Let us consider $\lambda \in (\lambda_{m+1}, \lambda_m)$. By (2.3) we have

$$\|\mathbf{y} - \hat{\boldsymbol{\mu}}_\lambda\|^2 = \mathbf{y}^T(\mathbf{I} - H_{\mathcal{B}_m})\mathbf{y} + \frac{\lambda^2}{4}\text{Sgn}_m^T(\mathbf{X}_{\mathcal{B}_m}^T\mathbf{X}_{\mathcal{B}_m})^{-1}\text{Sgn}_m, \tag{2.19}$$

where $H_{\mathcal{B}_m} = \mathbf{X}_{\mathcal{B}_m}(\mathbf{X}_{\mathcal{B}_m}^T\mathbf{X}_{\mathcal{B}_m})^{-1}\mathbf{X}_{\mathcal{B}_m}^T$. Thus we can conclude that $\|\mathbf{y} - \hat{\boldsymbol{\mu}}_\lambda\|^2$ is strictly increasing in the interval $(\lambda_{m+1}, \lambda_m)$. Moreover, the lasso estimates are continuous on $\lambda$, hence $\|\mathbf{y} - \hat{\boldsymbol{\mu}}_{\lambda_m}\|^2 > \|\mathbf{y} - \hat{\boldsymbol{\mu}}_\lambda\|^2 > \|\mathbf{y} - \hat{\boldsymbol{\mu}}_{\lambda_{m+1}}\|^2$. On the other hand, note that $\widehat{df}(\lambda) = |\mathcal{B}_m|\ \forall \lambda \in (\lambda_{m+1}, \lambda_m)$ and $|\mathcal{B}_m| \geq |\mathcal{B}(\lambda_{m+1})|$. Therefore the optimal choice of $\lambda$ in $[\lambda_{m+1}, \lambda_m)$ is $\lambda_{m+1}$, which means $\lambda(\text{optimal}) \in \{\lambda_m\}$.  □

According to Theorem 3, the optimal lasso model is immediately selected once we compute the entire lasso solution paths by the LARS algorithm. We can finish the whole fitting and tuning process with the computational cost of a single least squares fit.

**3. Efron's conjecture.** Efron et al. [5] first considered deriving the analytical form of the degrees of freedom of the lasso. They proposed a stagewise algorithm called LARS to compute the entire lasso solution paths. They also presented the following conjecture on the degrees of freedom of the lasso:

CONJECTURE 1. *Starting at step 0, let $m_k^{\text{last}}$ be the index of the last LARS-lasso sequence containing exactly $k$ nonzero predictors. Then $df(\hat{\boldsymbol{\mu}}_{m_k^{\text{last}}}) = k$.*

Note that Efron et al. [5] viewed the lasso as a forward stage-wise modeling algorithm and used the number of steps as the tuning parameter in the lasso: the lasso is regularized by early stopping. In the previous sections we regarded the lasso as a continuous penalization method with $\lambda$ as its regularization parameter. There is a subtle but important difference between the two views. The $\lambda$ value associated with $m_k^{\text{last}}$ is a random quantity. In the forward stage-wise modeling view of the lasso, the conjecture cannot be used for the degrees of freedom of the lasso at a general step $k$ for a prefixed $k$. This is simply because the number of LARS-lasso steps can exceed the number of all predictors (Efron et al. [5]). In contrast, the unbiasedness property of $\widehat{df}(\lambda)$ holds for all $\lambda$.

In this section we provide some justifications for the conjecture:

- We give a much more simplified proof than that in Efron et al. [5] to show that the conjecture is true under the positive cone condition.
- Our analysis also indicates that without the positive cone condition the conjecture can be wrong, although $k$ is a good approximation of $df(\hat{\boldsymbol{\mu}}_{m_k^{\text{last}}})$.
- We show that the conjecture works appropriately from the model selection perspective. If we use the conjecture to construct AIC (or BIC) to select the lasso fit, then the selected model is identical to that selected by AIC (or BIC) using the exact degrees of freedom results in Section 2.4.

First, we need to show that with probability one we can well define the last LARS-lasso sequence containing exactly $k$ nonzero predictors. Since the conjecture becomes a simple fact for the two trivial cases $k = 0$ and $k = p$, we only need to consider $k = 1, \ldots, p-1$. Let $\Lambda_k = \{m : |\mathcal{B}_{\lambda_m}| = k\}, k \in \{1, 2, \ldots, (p-1)\}$. Then $m_k^{\text{last}} = \sup(\Lambda_k)$. However, it may happen that for some $k$ there is no such $m$ with $|\mathcal{B}_{\lambda_m}| = k$. For example, if $\mathbf{y}$ is an equiangular vector of all $\{\mathbf{X}_j\}$, then the lasso estimates become the OLS estimates after just one step. So $\Lambda_k = \varnothing$ for $k = 2, \ldots, p-1$. The next lemma shows that the "*one at a time*" condition (Efron et al. [5]) holds almost everywhere; therefore $m_k^{\text{last}}$ is well defined almost surely.

LEMMA 7. *Let $\mathcal{W}_m(\mathbf{y})$ denote the set of predictors that are to be included in the active set at $\lambda_m$ and let $\mathcal{V}_m(\mathbf{y})$ be the set of predictors that are deleted from the active set at $\lambda_{m+1}$. Then $\exists$ a set $\widetilde{\mathcal{N}}_0$ which is a collection of finite many hyperplanes in $\mathbb{R}^n$. $\forall \mathbf{y} \in \mathbb{R}^n \setminus \widetilde{\mathcal{N}}_0$,*

(3.1) $\quad |\mathcal{W}_m(\mathbf{y})| \leq 1 \quad and \quad |\mathcal{V}_m(\mathbf{y})| \leq 1 \quad \forall m = 0, 1, \ldots, K(\mathbf{y})$.

$\mathbf{y} \in \mathbb{R}^n \setminus \widetilde{\mathcal{N}}_0$ is said to be a locally stable point for $\Lambda_k$, if $\forall \mathbf{y}'$ such that $\|\mathbf{y}' - \mathbf{y}\| \leq \varepsilon(\mathbf{y})$ for a small enough $\varepsilon(\mathbf{y})$, the effective set $\mathcal{B}(\lambda_{m_k^{\text{last}}})(\mathbf{y}') = \mathcal{B}(\lambda_{m_k^{\text{last}}})(\mathbf{y})$. Let $LS(k)$ be the set of all locally stable points.

The next lemma helps us evaluate $df(\hat{\boldsymbol{\mu}}_{m_k^{\text{last}}})$.

LEMMA 8. *Let $\hat{\boldsymbol{\mu}}_m(\mathbf{y})$ be the lasso fit at the transition point $\lambda_m$, $\lambda_m > 0$. Then for any $i \in \mathcal{W}_m$, we can write $\hat{\boldsymbol{\mu}}(m)$ as*

$$\hat{\boldsymbol{\mu}}_m(\mathbf{y}) = \left\{ \mathbf{H}_{\mathcal{B}(\lambda_m)} \right. \tag{3.2}$$
$$\left. - \frac{\mathbf{X}_{\mathcal{B}(\lambda_m)}^T (\mathbf{X}_{\mathcal{B}(\lambda_m)}^T \mathbf{X}_{\mathcal{B}(\lambda_m)}) \operatorname{Sgn}(\lambda_m) \mathbf{x}_i^T (\mathbf{I} - \mathbf{H}_{\mathcal{B}(\lambda_m)})}{\operatorname{Sgn}_i - \mathbf{x}_i^T \mathbf{X}_{\mathcal{B}(\lambda_m)}^T (\mathbf{X}_{\mathcal{B}(\lambda_m)}^T \mathbf{X}_{\mathcal{B}(\lambda_m)}) \operatorname{Sgn}(\lambda_m)} \right\} \mathbf{y}$$

$$=: \mathbf{S}_m(\mathbf{y})\mathbf{y}, \tag{3.3}$$

*where $\mathbf{H}_{\mathcal{B}(\lambda_m)}$ is the projection matrix on the subspace of $\mathbf{X}_{\mathcal{B}(\lambda_m)}$. Moreover*

$$\operatorname{tr}(\mathbf{S}_m(\mathbf{y})) = |\mathcal{B}(\lambda_m)|. \tag{3.4}$$

Note that $|\mathcal{B}(\lambda_{m_k^{\text{last}}})| = k$. Therefore, if $\mathbf{y} \in LS(k)$, then

$$\nabla \cdot \hat{\boldsymbol{\mu}}_{m_k^{\text{last}}}(\mathbf{y}) = \operatorname{tr}(\mathbf{S}_{m_k^{\text{last}}}(\mathbf{y})) = k. \tag{3.5}$$

If the positive cone condition holds then the lasso solution paths are monotone (Efron et al. [5]), hence Lemma 7 implies that $LS(k)$ is a set of full measure. Then by Lemma 8 we know that $df(m_k^{\text{last}}) = k$. However, it should be pointed out that $k - df(m_k^{\text{last}})$ can be nonzero for some $k$ when the positive cone condition is violated. Here we present an explicit example to show this point. We consider the synthetic model in Section 2.3. Note that the positive cone condition is violated on the 64 predictors [5]. As done in Section 2.3, the exact value of $df(m_k^{\text{last}})$ can be computed by Monte Carlo and then we evaluate the bias $k - df(m_k^{\text{last}})$. In the synthetic model (2.12) the signal/noise ratio $\frac{\operatorname{Var}(\mathbf{X}\hat{\beta}_{\text{ols}})}{\hat{\sigma}_{\text{ols}}^2}$ is about 1.25. We repeated the same simulation procedure with $(\beta = \hat{\beta}_{\text{ols}}, \sigma = \frac{\hat{\sigma}_{\text{ols}}}{10})$ in the synthetic model and the corresponding signal/noise ratio became 125. As shown clearly in Figure 4, the bias $k - df(m_k^{\text{last}})$ is not zero for some $k$. However, even if the bias exists, its maximum magnitude is less than one, regardless of the size of the signal/noise ratio, which suggests that $k$ is a good estimate of $df(m_k^{\text{last}})$.

Let us pretend the conjecture is true in all situations and then define the model selection criteria as

$$\frac{\|\mathbf{y} - \hat{\boldsymbol{\mu}}_{m_k^{\text{last}}}\|^2}{n\sigma^2} + \frac{w_n}{n} k. \tag{3.6}$$

$w_n = 2$ for AIC and $w_n = \log(n)$ for BIC. Treat $k$ as the tuning parameter of the lasso. We need to find $k(\text{optimal})$ such that

$$k(\text{optimal}) = \arg\min_k \frac{\|\mathbf{y} - \hat{\boldsymbol{\mu}}_{m_k^{\text{last}}}\|^2}{n\sigma^2} + \frac{w_n}{n} k. \tag{3.7}$$

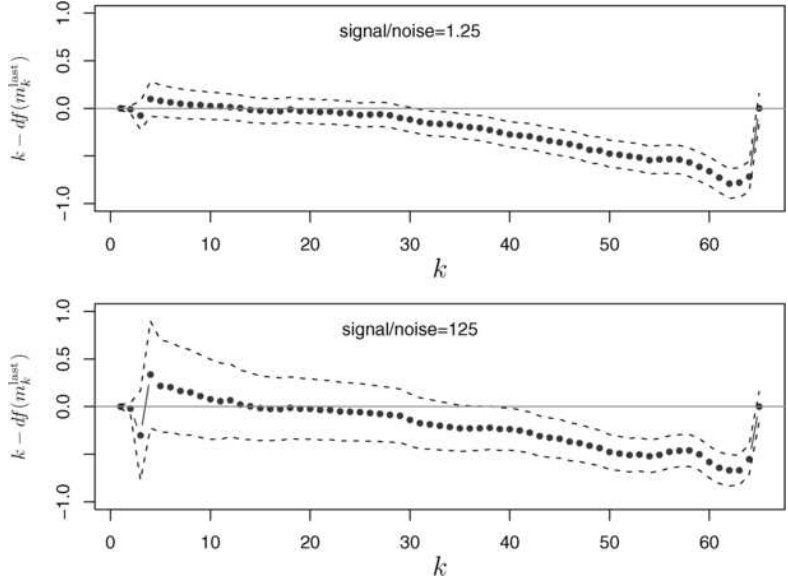

FIG. 4. $B = 20000$ replications were used to assess the bias of $\widehat{df}(m_k^{\text{last}}) = k$. The 95% point-wise confidence intervals are indicated by the thin dashed lines. This simulation suggests that when the positive cone condition is violated, $df(m_k^{\text{last}}) \neq k$ for some $k$. However, the bias is small (the maximum absolute bias is about 0.8), regardless of the size of the signal/noise ratio.

Suppose $\lambda^* = \lambda(\text{optimal})$ and $k^* = k(\text{optimal})$. Theorem 3 implies that the models selected by (2.17) and (3.7) coincide, that is, $\hat{\boldsymbol{\mu}}_{\lambda^*} = \hat{\boldsymbol{\mu}}_{m_{k^*}^{\text{last}}}$. This observation suggests that although the conjecture is not always true, it actually works appropriately for the purpose of model selection.

**4. Proofs of the lemmas.** First, let us introduce the following matrix representation of the divergence. Let $\frac{\partial \hat{\boldsymbol{\mu}}}{\partial \mathbf{y}}$ be a $n \times n$ matrix whose elements are

$$\left(\frac{\partial \hat{\boldsymbol{\mu}}}{\partial \mathbf{y}}\right)_{i,j} = \frac{\partial \hat{\mu}_i}{\partial y_j}, \qquad i, j = 1, 2, \ldots, n. \tag{4.1}$$

Then we can write

$$\nabla \cdot \hat{\boldsymbol{\mu}} = \text{tr}\left(\frac{\partial \hat{\boldsymbol{\mu}}}{\partial \mathbf{y}}\right). \tag{4.2}$$

The above trace expression will be used repeatedly.

PROOF OF LEMMA 1. Let

$$\ell(\beta, \mathbf{y}) = \left\| \mathbf{y} - \sum_{j=1}^{p} \mathbf{x}_j \beta_j \right\|^2 + \lambda \sum_{j=1}^{p} |\beta_j|. \tag{4.3}$$

Given $\mathbf{y}$, $\hat{\beta}(\lambda)$ is the minimizer of $\ell(\beta, \mathbf{y})$. For those $j \in \mathcal{B}_m$ we must have $\frac{\partial \ell(\beta, \mathbf{y})}{\partial \beta_j} = 0$, that is,

$$(4.4) \quad -2\mathbf{x}_j^T\left(\mathbf{y} - \sum_{j=1}^p \mathbf{x}_j \hat{\beta}(\lambda)_j\right) + \lambda \operatorname{Sgn}(\hat{\beta}(\lambda)_j) = 0, \qquad \text{for } j \in \mathcal{B}_m.$$

Since $\hat{\beta}(\lambda)_i = 0$ for all $i \notin \mathcal{B}_m$, then $\sum_{j=1}^p \mathbf{x}_j \hat{\beta}(\lambda)_j = \sum_{j \in \mathcal{B}_\lambda} \mathbf{x}_j \hat{\beta}(\lambda)_j$. Thus the equations in (4.4) become

$$(4.5) \quad -2\mathbf{X}_{\mathcal{B}_m}^T(\mathbf{y} - \mathbf{X}_{\mathcal{B}_m} \hat{\beta}(\lambda)_{\mathcal{B}_m}) + \lambda \operatorname{Sgn}_m = 0,$$

which gives (2.3). $\square$

PROOF OF LEMMA 2. We adopt the matrix notation used in SPLUS: $\mathbf{M}[i, \cdot]$ means the $i$th row of $\mathbf{M}$. $i_{\mathrm{add}}$ joins $\mathcal{B}_m$ at $\lambda_m$; then $\hat{\beta}(\lambda_m)_{i_{\mathrm{add}}} = 0$. Consider $\hat{\beta}(\lambda)$ for $\lambda \in (\lambda_{m+1}, \lambda_m)$. Lemma 1 gives

$$(4.6) \quad \hat{\beta}(\lambda)_{\mathcal{B}_m} = (\mathbf{X}_{\mathcal{B}_m}^T \mathbf{X}_{\mathcal{B}_m})^{-1}\left(\mathbf{X}_{\mathcal{B}_m}^T \mathbf{y} - \frac{\lambda}{2} \operatorname{Sgn}_m\right).$$

By the continuity of $\hat{\beta}(\lambda)_{i_{\mathrm{add}}}$, taking the limit of the $i^*$th element of (4.6) as $\lambda \to \lambda_m - 0$, we have

$$(4.7) \quad 2\{(\mathbf{X}_{\mathcal{B}_m}^T \mathbf{X}_{\mathcal{B}_m})^{-1}[i^*, \cdot] \mathbf{X}_{\mathcal{B}_m}^T\} \mathbf{y} = \lambda_m \{(\mathbf{X}_{\mathcal{B}_m}^T \mathbf{X}_{\mathcal{B}_m})^{-1}[i^*, \cdot] \operatorname{Sgn}_m\}.$$

The second $\{\cdot\}$ is a nonzero scalar, otherwise $\hat{\beta}(\lambda)_{i_{\mathrm{add}}} = 0$ for all $\lambda \in (\lambda_{m+1}, \lambda_m)$, which contradicts the assumption that $i_{\mathrm{add}}$ becomes a member of the active set $\mathcal{B}_m$. Thus we have

$$(4.8) \quad \lambda_m = \left\{2\frac{(\mathbf{X}_{\mathcal{B}_m}^T \mathbf{X}_{\mathcal{B}_m})^{-1}[i^*, \cdot]}{(\mathbf{X}_{\mathcal{B}_m}^T \mathbf{X}_{\mathcal{B}_m})^{-1}[i^*, \cdot] \operatorname{Sgn}_m}\right\} \mathbf{X}_{\mathcal{B}_m}^T \mathbf{y} =: v(\mathcal{B}_m, i^*) \mathbf{X}_{\mathcal{B}_m}^T \mathbf{y},$$

where $v(\mathcal{B}_m, i^*) = \{2((\mathbf{X}_{\mathcal{B}_m}^T \mathbf{X}_{\mathcal{B}_m})^{-1}[i^*, \cdot])/((\mathbf{X}_{\mathcal{B}_m}^T \mathbf{X}_{\mathcal{B}_m})^{-1}[i^*, \cdot] \operatorname{Sgn}_m)\}$. Rearranging (4.8), we get (2.4).

Similarly, if $j_{\mathrm{drop}}$ is a dropped index at $\lambda_{m+1}$, we take the limit of the $j^*$th element of (4.6) as $\lambda \to \lambda_{m+1} + 0$ to conclude that

$$(4.9) \quad \lambda_{m+1} = \left\{2\frac{(\mathbf{X}_{\mathcal{B}_m}^T \mathbf{X}_{\mathcal{B}_m})^{-1}[j^*, \cdot]}{(\mathbf{X}_{\mathcal{B}_m}^T \mathbf{X}_{\mathcal{B}_m})^{-1}[j^*, \cdot] \operatorname{Sgn}_m}\right\} \mathbf{X}_{\mathcal{B}_m}^T \mathbf{y} =: v(\mathcal{B}_m, j^*) \mathbf{X}_{\mathcal{B}_m}^T \mathbf{y},$$

where $v(\mathcal{B}_m, j^*) = \{2((\mathbf{X}_{\mathcal{B}_m}^T \mathbf{X}_{\mathcal{B}_m})^{-1}[j^*, \cdot])/((\mathbf{X}_{\mathcal{B}_m}^T \mathbf{X}_{\mathcal{B}_m})^{-1}[j^*, \cdot] \operatorname{Sgn}_m)\}$. Rearranging (4.9), we get (2.5). $\square$

PROOF OF LEMMA 3. Suppose for some $\mathbf{y}$ and $m$, $\lambda = \lambda(\mathbf{y})_m$. $\lambda > 0$ means $m$ is not the last lasso step. By Lemma 2 we have

$$(4.10) \quad \lambda = \lambda_m = \{v(\mathcal{B}_m, i^*) \mathbf{X}_{\mathcal{B}_m}^T\} \mathbf{y} =: \alpha(\mathcal{B}_m, i^*) \mathbf{y}.$$

Obviously $\alpha(\mathcal{B}_m, i^*) = v(\mathcal{B}_m, i^*)\mathbf{X}_{\mathcal{B}_m}^T$ is a nonzero vector. Now let $\alpha_\lambda$ be the totality of $\alpha(\mathcal{B}_m, i^*)$ by considering all the possible combinations of $\mathcal{B}_m$, $i^*$ and the sign vector $\text{Sgn}_m$. $\alpha_\lambda$ depends only on $\mathbf{X}$ and is a finite set, since at most $p$ predictors are available. Thus $\forall \alpha \in \alpha_\lambda$, $\alpha \mathbf{y} = \lambda$ defines a hyperplane in $\mathbb{R}^n$. We define

$$\mathcal{N}_\lambda = \{\mathbf{y} : \alpha \mathbf{y} = \lambda \text{ for some } \alpha \in \alpha_\lambda\} \quad \text{and} \quad \mathcal{G}_\lambda = \mathbb{R}^n \setminus \mathcal{N}_\lambda.$$

Then on $\mathcal{G}_\lambda$ (4.10) is impossible. $\square$

PROOF OF LEMMA 4. For writing convenience we omit the subscript $\lambda$. Let $\hat{\beta}(\mathbf{y})_{\text{ols}} = (\mathbf{X}^T \mathbf{X})^{-1} \mathbf{X}^T \mathbf{y}$ be the OLS estimates. Note that we always have the inequality

$$|\hat{\beta}(\mathbf{y})|_1 \leq |\hat{\beta}(\mathbf{y})_{\text{ols}}|_1. \tag{4.11}$$

Fix an arbitrary $\mathbf{y}_0$ and consider a sequence of $\{\mathbf{y}_n\}$ ($n = 1, 2, \ldots$) such that $\mathbf{y}_n \to \mathbf{y}_0$. Since $\mathbf{y}_n \to \mathbf{y}_0$, we can find a $Y$ such that $\|\mathbf{y}_n\| \leq Y$ for all $n = 0, 1, 2, \ldots$. Consequently $\|\hat{\beta}(\mathbf{y}_n)_{\text{ols}}\| \leq B$ for some upper bound $B$ ($B$ is determined by $\mathbf{X}$ and $Y$). By Cauchy's inequality and (4.11), we have $|\hat{\beta}(\mathbf{y}_n)|_1 \leq \sqrt{pB}$ for all $n = 0, 1, 2, \ldots$. Thus to show $\hat{\beta}(\mathbf{y}_n) \to \hat{\beta}(\mathbf{y}_0)$, it is equivalent to show that for every converging subsequence of $\{\hat{\beta}(\mathbf{y}_n)\}$, say $\{\hat{\beta}(\mathbf{y}_{n_k})\}$, the subsequence converges to $\hat{\beta}(\mathbf{y})$. Now suppose $\hat{\beta}(\mathbf{y}_{n_k})$ converges to $\hat{\beta}_\infty$ as $n_k \to \infty$. We show $\hat{\beta}_\infty = \hat{\beta}(\mathbf{y}_0)$. The lasso criterion $\ell(\beta, \mathbf{y})$ is written in (4.3). Let $\Delta \ell(\beta, \mathbf{y}, \mathbf{y}') = \ell(\beta, \mathbf{y}) - \ell(\beta, \mathbf{y}')$. By the definition of $\hat{\beta}_{n_k}$, we must have

$$\ell(\hat{\beta}(\mathbf{y}_0), \mathbf{y}_{n_k}) \geq \ell(\hat{\beta}(\mathbf{y}_{n_k}), \mathbf{y}_{n_k}). \tag{4.12}$$

Then (4.12) gives

$$\begin{aligned}
\ell(\hat{\beta}(\mathbf{y}_0), \mathbf{y}_0) &= \ell(\hat{\beta}(\mathbf{y}_0), \mathbf{y}_{n_k}) + \Delta \ell(\hat{\beta}(\mathbf{y}_0), \mathbf{y}_0, \mathbf{y}_{n_k}) \\
&\geq \ell(\hat{\beta}(\mathbf{y}_{n_k}), \mathbf{y}_{n_k}) + \Delta \ell(\hat{\beta}(\mathbf{y}_0), \mathbf{y}_0, \mathbf{y}_{n_k}) \\
&= \ell(\hat{\beta}(\mathbf{y}_{n_k}), \mathbf{y}_0) + \Delta \ell(\hat{\beta}(\mathbf{y}_{n_k}), \mathbf{y}_{n_k}, \mathbf{y}_0) \\
&\quad + \Delta \ell(\hat{\beta}(\mathbf{y}_0), \mathbf{y}_0, \mathbf{y}_{n_k}).
\end{aligned} \tag{4.13}$$

We observe

$$\begin{aligned}
\Delta \ell(\hat{\beta}(\mathbf{y}_{n_k}), \mathbf{y}_{n_k}, \mathbf{y}_0) &+ \Delta \ell(\hat{\beta}(\mathbf{y}_0), \mathbf{y}_0, \mathbf{y}_{n_k}) \\
&= 2(\mathbf{y}_0 - \mathbf{y}_{n_k}) \mathbf{X}^T (\hat{\beta}(\mathbf{y}_{n_k}) - \hat{\beta}(\mathbf{y}_0)).
\end{aligned} \tag{4.14}$$

Let $n_k \to \infty$; the right-hand side of (4.14) goes to zero. Moreover, $\ell(\hat{\beta}(\mathbf{y}_{n_k}), \mathbf{y}_0) \to \ell(\hat{\beta}_\infty, \mathbf{y}_0)$. Therefore (4.13) reduces to

$$\ell(\hat{\beta}(\mathbf{y}_0), \mathbf{y}_0) \geq \ell(\hat{\beta}_\infty, \mathbf{y}_0).$$

However, $\hat{\beta}(\mathbf{y}_0)$ is the unique minimizer of $\ell(\beta, \mathbf{y}_0)$, and thus $\hat{\beta}_\infty = \hat{\beta}(\mathbf{y}_0)$.
$\square$

PROOF OF LEMMA 5. Fix an arbitrary $\mathbf{y}_0 \in \mathcal{G}_\lambda$. Denote by $\text{Ball}(\mathbf{y}, r)$ the $n$-dimensional ball with center $\mathbf{y}$ and radius $r$. Note that $\mathcal{G}_\lambda$ is an open set, so we can choose a small enough $\varepsilon$ such that $\text{Ball}(\mathbf{y}_0, \varepsilon) \subset \mathcal{G}_\lambda$. Fix $\varepsilon$. Suppose $\mathbf{y}_n \to \mathbf{y}$ as $n \to \infty$. Then without loss of generality we can assume $\mathbf{y}_n \in \text{Ball}(\mathbf{y}_0, \varepsilon)$ for all $n$. So $\lambda$ is not a transition point for any $\mathbf{y}_n$.

By definition $\hat{\beta}(\mathbf{y}_0)_j \neq 0$ for all $j \in \mathcal{B}(\mathbf{y}_0)$. Then Lemma 4 says that $\exists$ an $N_1$, and as long as $n > N_1$, we have $\hat{\beta}(\mathbf{y}_n)_j \neq 0$ and $\text{Sgn}(\hat{\beta}(\mathbf{y}_n)) = \text{Sgn}(\hat{\beta}(\mathbf{y}_n))$, for all $j \in \mathcal{B}(\mathbf{y}_0)$. Thus $\mathcal{B}(\mathbf{y}_0) \subseteq \mathcal{B}(\mathbf{y}_n) \ \forall n > N_1$.

On the other hand, we have the *equiangular* conditions (Efron et al. [5])

$$\lambda = 2|\mathbf{x}_j^T(\mathbf{y}_0 - \mathbf{X}\hat{\beta}(\mathbf{y}_0))| \qquad \forall j \in \mathcal{B}(\mathbf{y}_0), \tag{4.15}$$

$$\lambda > 2|\mathbf{x}_j^T(\mathbf{y}_0 - \mathbf{X}\hat{\beta}(\mathbf{y}_0))| \qquad \forall j \notin \mathcal{B}(\mathbf{y}_0). \tag{4.16}$$

Using Lemma 4 again, we conclude that $\exists$ an $N > N_1$ such that $\forall j \notin \mathcal{B}(\mathbf{y}_0)$ the strict inequalities (4.16) hold for $\mathbf{y}_n$ provided $n > N$. Thus $\mathcal{B}^c(\mathbf{y}_0) \subseteq \mathcal{B}^c(\mathbf{y}_n) \ \forall n > N$. Therefore we have $\mathcal{B}(\mathbf{y}_n) = \mathcal{B}(\mathbf{y}_0) \ \forall n > N$. Then the local constancy of the sign vector follows the continuity of $\hat{\beta}(\mathbf{y})$. $\square$

PROOF OF LEMMA 6. If $\lambda = 0$, then the lasso fit is just the OLS fit. The conclusions are easy to verify. So we focus on $\lambda > 0$. Fix an $\mathbf{y}$. Choose a small enough $\varepsilon$ such that $\text{Ball}(\mathbf{y}, \varepsilon) \subset \mathcal{G}_\lambda$.

Since $\lambda$ is not any transition point, using (2.3) we observe

$$\hat{\boldsymbol{\mu}}_\lambda(\mathbf{y}) = \mathbf{X}\hat{\beta}(\mathbf{y}) = \mathbf{H}_\lambda(\mathbf{y})\mathbf{y} - \lambda\boldsymbol{\omega}_\lambda(\mathbf{y}), \tag{4.17}$$

where $\mathbf{H}_\lambda(\mathbf{y}) = \mathbf{X}_{\mathcal{B}_\lambda}(\mathbf{X}_{\mathcal{B}_\lambda}^T \mathbf{X}_{\mathcal{B}_\lambda})^{-1}\mathbf{X}_{\mathcal{B}_\lambda}^T$ is the projection matrix on the space $\mathbf{X}_{\mathcal{B}_\lambda}$ and $\boldsymbol{\omega}_\lambda(\mathbf{y}) = \frac{1}{2}\mathbf{X}_{\mathcal{B}_\lambda}(\mathbf{X}_{\mathcal{B}_\lambda}^T \mathbf{X}_{\mathcal{B}_\lambda})^{-1}\text{Sgn}_{\mathcal{B}_\lambda}$. Consider $\|\Delta\mathbf{y}\| < \varepsilon$. Similarly, we get

$$\hat{\boldsymbol{\mu}}_\lambda(\mathbf{y} + \Delta\mathbf{y}) = \mathbf{H}_\lambda(\mathbf{y} + \Delta\mathbf{y})(\mathbf{y} + \Delta\mathbf{y}) - \lambda\boldsymbol{\omega}_\lambda(\mathbf{y} + \Delta\mathbf{y}). \tag{4.18}$$

Lemma 5 says that we can further let $\varepsilon$ be sufficiently small such that both the effective set $\mathcal{B}_\lambda$ and the sign vector $\text{Sgn}_\lambda$ stay constant in $\text{Ball}(\mathbf{y}, \varepsilon)$. Now fix $\varepsilon$. Hence if $\|\Delta\mathbf{y}\| < \varepsilon$, then

$$\mathbf{H}_\lambda(\mathbf{y} + \Delta\mathbf{y}) = \mathbf{H}_\lambda(\mathbf{y}) \quad \text{and} \quad \boldsymbol{\omega}_\lambda(\mathbf{y} + \Delta\mathbf{y}) = \boldsymbol{\omega}_\lambda(\mathbf{y}). \tag{4.19}$$

Then (4.17) and (4.18) give

$$\hat{\boldsymbol{\mu}}_\lambda(\mathbf{y} + \Delta\mathbf{y}) - \hat{\boldsymbol{\mu}}_\lambda(\mathbf{y}) = \mathbf{H}_\lambda(\mathbf{y})\Delta\mathbf{y}. \tag{4.20}$$

But since $\|\mathbf{H}_\lambda(\mathbf{y})\Delta\mathbf{y}\| \leq \|\Delta\mathbf{y}\|$, (2.6) is proved.

By the local constancy of $H(\mathbf{y})$ and $\omega(\mathbf{y})$, we have

$$\frac{\partial \hat{\boldsymbol{\mu}}_\lambda(\mathbf{y})}{\partial \mathbf{y}} = \mathbf{H}_\lambda(\mathbf{y}). \tag{4.21}$$

Then the trace formula (4.2) implies that

$$\nabla \cdot \hat{\boldsymbol{\mu}}_\lambda(\mathbf{y}) = \mathrm{tr}(\mathbf{H}_\lambda(\mathbf{y})) = |\mathcal{B}_\lambda|. \tag{4.22} \qquad \square$$

PROOF OF LEMMA 7. Suppose at step $m$, $|\mathcal{W}_m(\mathbf{y})| \geq 2$. Let $i_{\mathrm{add}}$ and $j_{\mathrm{add}}$ be two of the predictors in $\mathcal{W}_m(\mathbf{y})$, and let $i^*_{\mathrm{add}}$ and $j^*_{\mathrm{add}}$ be their indices in the current active set $\mathcal{A}$. Note the current active set $\mathcal{A}$ is $\mathcal{B}_m$ in Lemma 2. Hence we have

$$\lambda_m = v[\mathcal{A}, i^*] \mathbf{X}_\mathcal{A}^T \mathbf{y} \quad \text{and} \quad \lambda_m = v[\mathcal{A}, j^*] \mathbf{X}_\mathcal{A}^T \mathbf{y}. \tag{4.23}$$

Therefore

$$0 = \{[v(\mathcal{A}, i^*_{\mathrm{add}}) - v(\mathcal{A}, j^*_{\mathrm{add}})] \mathbf{X}_\mathcal{A}^T\} \mathbf{y} =: \alpha_{\mathrm{add}} \mathbf{y}. \tag{4.24}$$

We claim $\alpha_{\mathrm{add}} = [v(\mathcal{A}, i^*_{\mathrm{add}}) - v(\mathcal{A}, j^*_{\mathrm{add}})] \mathbf{X}_\mathcal{A}^T$ is not a zero vector. Otherwise, since $\{\mathbf{X}_j\}$ are linearly independent, $\alpha_{\mathrm{add}} = 0$ forces $v(\mathcal{A}, i^*_{\mathrm{add}}) - v(\mathcal{A}, j^*_{\mathrm{add}}) = 0$. Then we have

$$\frac{(\mathbf{X}_\mathcal{A}^T \mathbf{X}_\mathcal{A})^{-1}[i^*, \cdot]}{(\mathbf{X}_\mathcal{A}^T \mathbf{X}_\mathcal{A})^{-1}[i^*, \cdot] \mathrm{Sgn}_\mathcal{A}} = \frac{(\mathbf{X}_\mathcal{A}^T \mathbf{X}_\mathcal{A})^{-1}[j^*, \cdot]}{(\mathbf{X}_\mathcal{A}^T \mathbf{X}_\mathcal{A})^{-1}[i^*, \cdot] \mathrm{Sgn}_\mathcal{A}}, \tag{4.25}$$

which contradicts the fact $(\mathbf{X}_\mathcal{A}^T \mathbf{X}_\mathcal{A})^{-1}$ is a full rank matrix.

Similarly, if $i_{\mathrm{drop}}$ and $j_{\mathrm{drop}}$ are dropped predictors, then

$$0 = \{[v(\mathcal{A}, i^*_{\mathrm{drop}}) - v(\mathcal{A}, j^*_{\mathrm{drop}})] \mathbf{X}_\mathcal{A}^T\} \mathbf{y} =: \alpha_{\mathrm{drop}} \mathbf{y}, \tag{4.26}$$

and $\alpha_{\mathrm{drop}} = [v(\mathcal{A}, i^*_{\mathrm{drop}}) - v(\mathcal{A}, j^*_{\mathrm{drop}})] \mathbf{X}_\mathcal{A}^T$ is a nonzero vector.

Let $M_0$ be the totality of $\alpha_{\mathrm{add}}$ and $\alpha_{\mathrm{drop}}$ by considering all the possible combinations of $\mathcal{A}$, $(i_{\mathrm{add}}, j_{\mathrm{add}})$, $(i_{\mathrm{drop}}, j_{\mathrm{drop}})$ and $\mathrm{Sgn}_\mathcal{A}$. Clearly $M_0$ is a finite set and depends only on $\mathbf{X}$. Let

$$\widetilde{\mathcal{N}_0} = \{y : \alpha y = 0 \text{ for some } \alpha \in M_0\}. \tag{4.27}$$

Then on $\mathbb{R}^n \setminus \widetilde{\mathcal{N}_0}$ the conclusion holds. $\square$

PROOF OF LEMMA 8. Note that $\hat{\beta}(\lambda)$ is continuous on $\lambda$. Using (4.4) in Lemma 1 and taking the limit of $\lambda \to \lambda_m$, we have

$$-2\mathbf{x}_j^T \left( \mathbf{y} - \sum_{j=1}^p \mathbf{x}_j \hat{\beta}(\lambda_m)_j \right) + \lambda_m \mathrm{Sgn}(\hat{\beta}(\lambda_m)_j) = 0, \qquad \text{for } j \in \mathcal{B}(\lambda_m). \tag{4.28}$$

However, $\sum_{j=1}^p \mathbf{x}_j \hat{\beta}(\lambda_m)_j = \sum_{j \in \mathcal{B}(\lambda_m)} \mathbf{x}_j \hat{\beta}(\lambda_m)_j$. Thus we have

$$\hat{\beta}(\lambda_m) = (\mathbf{X}_{\mathcal{B}(\lambda_m)}^T \mathbf{X}_{\mathcal{B}(\lambda_m)})^{-1} \left( \mathbf{X}_{\mathcal{B}(\lambda_m)}^T \mathbf{y} - \frac{\lambda_m}{2} \mathrm{Sgn}(\lambda_m) \right). \tag{4.29}$$

Hence

$$\hat{\boldsymbol{\mu}}_m(\mathbf{y}) = \mathbf{X}_{\mathcal{B}(\lambda_m)}(\mathbf{X}_{\mathcal{B}(\lambda_m)}^T \mathbf{X}_{\mathcal{B}(\lambda_m)})^{-1}\left(\mathbf{X}_{\mathcal{B}(\lambda_m)}^T \mathbf{y} - \frac{\lambda_m}{2}\operatorname{Sgn}(\lambda_m)\right)$$
(4.30)
$$= \mathbf{H}_{\mathcal{B}(\lambda_m)}\mathbf{y} - \mathbf{X}_{\mathcal{B}(\lambda_m)}(\mathbf{X}_{\mathcal{B}(\lambda_m)}^T \mathbf{X}_{\mathcal{B}(\lambda_m)})^{-1}\operatorname{Sgn}(\lambda_m)\frac{\lambda_m}{2}.$$

Since $i \in \mathcal{W}_m$, we must have the *equiangular* condition

$$\operatorname{Sgn}_i \mathbf{x}_i^T(\mathbf{y} - \hat{\boldsymbol{\mu}}(m)) = \frac{\lambda_m}{2}. \tag{4.31}$$

Substituting (4.30) into (4.31), we solve $\lambda_m/2$ and obtain

$$\frac{\lambda_m}{2} = \frac{\mathbf{x}_i^T(\mathbf{I} - \mathbf{H}_{\mathcal{B}(\lambda_m)})\mathbf{y}}{\operatorname{Sgn}_i - \mathbf{x}_i^T \mathbf{X}_{\mathcal{B}(\lambda_m)}^T(\mathbf{X}_{\mathcal{B}(\lambda_m)}^T \mathbf{X}_{\mathcal{B}(\lambda_m)})\operatorname{Sgn}(\lambda_m)}. \tag{4.32}$$

Then putting (4.32) back to (4.30) yields (3.2).

Using the identity $\operatorname{tr}(AB) = \operatorname{tr}(BA)$, we observe

$$\operatorname{tr}(\mathbf{S}_m(\mathbf{y}) - \mathbf{H}_{\mathcal{B}(\lambda_m)}) = \operatorname{tr}\left(\frac{(\mathbf{X}_{\mathcal{B}(\lambda_m)}^T \mathbf{X}_{\mathcal{B}(\lambda_m)})\operatorname{Sgn}(\lambda_m)\mathbf{x}_i^T(\mathbf{I} - H_{\mathcal{B}(\lambda_m)})\mathbf{X}_{\mathcal{B}(\lambda_m)}^T}{\operatorname{Sgn}_i - \mathbf{x}_i^T \mathbf{X}_{\mathcal{B}(\lambda_m)}^T(\mathbf{X}_{\mathcal{B}(\lambda_m)}^T \mathbf{X}_{\mathcal{B}(\lambda_m)})\operatorname{Sgn}(\lambda_m)}\right)$$
$$= \operatorname{tr}(0) = 0.$$

So $\operatorname{tr}(\mathbf{S}_m(\mathbf{y})) = \operatorname{tr}(\mathbf{H}_{\mathcal{B}(\lambda_m)}) = |\mathcal{B}(\lambda_m)|.$ □

**5. Discussion.** In this article we have proven that the number of nonzero coefficients is an unbiased estimate of the degrees of freedom of the lasso. The unbiased estimator is also consistent. We think it is a neat yet surprising result. Even in other sparse modeling methods, there is no such clean relationship between the number of nonzero coefficients and the degrees of freedom. For example, the number of nonzero coefficients is not an unbiased estimate of the degrees of freedom of the elastic net (Zou [26]). Another possible counterexample is the SCAD (Fan and Li [6]) whose solution is even more complex than the lasso. Note that with orthogonal predictors, the SCAD estimates can be obtained by the SCAD shrinkage formula (Fan and Li [6]). Then it is not hard to check that with orthogonal predictors the number of nonzero coefficients in the SCAD estimates cannot be an unbiased estimate of its degrees of freedom.

The techniques developed in this article can be applied to derive the degrees of freedom of other nonlinear estimating procedures, especially when the estimates have piece-wise linear solution paths. Gunter and Zhu [9] used our arguments to derive an unbiased estimate of the degrees of freedom of support vector regression. Zhao, Rocha and Yu [25] derived an unbiased estimate of the degrees of freedom of the regularized estimates using the CAP penalties.

Bühlmann and Yu [2] defined the degrees of freedom of $L_2$ boosting as the trace of the product of a series of linear smoothers. Their approach takes advantage of the closed-form expression for the $L_2$ fit at each boosting stage. It is now well known that $\varepsilon$-$L_2$ boosting is (almost) identical to the lasso (Hastie, Tibshirani and Friedman [11], Efron et al. [5]). Their work provides another look at the degrees of freedom of the lasso. However, it is not clear whether their definition agrees with the SURE definition. This could be another interesting topic for future research.

**Acknowledgments.** Hui Zou sincerely thanks Brad Efron, Yuhong Yang and Xiaotong Shen for their encouragement and suggestions. We sincerely thank the Co-Editor Jianqing Fan, an Associate Editor and two referees for helpful comments which greatly improved the manuscript.

H. ZOU
SCHOOL OF STATISTICS
UNIVERSITY OF MINNESOTA
MINNEAPOLIS, MINNESOTA 55455
USA
E-MAIL: hzou@stat.umn.edu

T. HASTIE
R. TIBSHIRANI
DEPARTMENT OF STATISTICS
STANFORD UNIVERSITY
STANFORD, CALIFORNIA 94305
USA
E-MAIL: hastie@stanford.edu
       tibs@stanford.edu